\documentclass{cjour} 
\usepackage{amsmath,amssymb,array,delarray,latexsym,wrapfig,epsfig} 

\newlength{\sh}
\newlength{\jmr}
\newlength{\jfc}
\newlength{\bernd}
\newlength{\jones}
\newlength{\mati}
\newlength{\tung}
\newlength{\sil}
\newlength{\koi}
\newlength{\gala}
\newtheorem{cor}{Corollary}
\newtheorem{prop}{Proposition}
\newtheorem{dfn}{Definition}

\newtheorem{main}{Theorem} 
\newtheorem{thm}[main]{Theorem}
\newtheorem{rem}{Remark}	
\newtheorem{ex}{Example}

\renewcommand{\mod}{\mathbf{mod}} 
\newcommand{\pspa}{{$\mathbf{PSPACE}$}} 
 
\newcommand{\am}{{$\mathbf{AM}$}} 

\newcommand{\np}{{\mathbf{NP}}}
\newcommand{\conp}{{\mathbf{coNP}}}
\newcommand{\corp}{{\mathbf{coRP}}}

\newcommand{\pp}{\mathbf{P}}
\newcommand{\nc}{\mathbf{NC}}
\newcommand{\expt}{{$\mathbf{EXPTIME}$}}

\newcommand{\cA}{\mathcal{A}}
\newcommand{\cD}{{\Delta}}

\newcommand{\cO}{\mathcal{O}}

\newcommand{\size}{\mathrm{size}}
\newcommand{\thth}{{\underline{\mathrm{th}}}}
\newcommand{\rd}{ {\underline{ \mathrm{rd} } }  }
\newcommand{\st}{ {\underline{ \mathrm{st} } }  }

\newcommand{\Pro}{\mathbb{P}}

\newcommand{\Q}{\mathbb{Q}}
\newcommand{\R}{\mathbb{R}}
\newcommand{\C}{\mathbb{C}}
\newcommand{\N}{\mathbb{N}}
\newcommand{\Z}{\mathbb{Z}}

\newcommand{\newt}{\mathrm{Newt}}

\newcommand{\Qn}{\Q^n} 
\newcommand{\Rn}{\R^n}

\newcommand{\Cn}{\C^n}

\renewcommand{\qed}{$\blacksquare$}

\newcommand{\cR}{\mathcal{R}}

\newcommand{\bO}{\mathbf{O}}
\newcommand{\vol}{\mathrm{Vol}}
\newcommand{\cd}{{\mathfrak{d}}}
\newcommand{\lie}{\mathrm{Li}}

\begin{document}

\authorrunninghead{J.\ Maurice Rojas} 
\titlerunninghead{Computational Arithmetic Geometry I} 
\title{Computational Arithmetic Geometry I: Sentences Nearly in the Polynomial 
Hierarchy\thanks{ 
An extended abstract of this work appeared earlier in the Proceedings of the 
31$^\st$ Annual ACM Symposium on Theory of Computing (STOC, May 1--4, 1999, 
Atlanta, Georgia), 527--536, ACM Press, 1999. This research was 
supported by Hong Kong UGC Grant \#9040469-730. }
} 

\author{J.\ Maurice Rojas\thanks{URL: {\tt 
http://math.cityu.edu.hk/\~{}mamrojas} } }  
\email{mamrojas@math.cityu.edu.hk}

\affil{Department of Mathematics, Texas A\&M University, College Station,
Texas \ 77843-3368, USA (after January 2001).\thanks{
Department of Mathematics, City University of Hong Kong,
83 Tat Chee Avenue, Kowloon, HONG KONG (until December 2000).}}

\dedication{Dedicated to Gretchen Davis.}  

\date{August 25, 2000} 

\abstract{ 
We consider the average-case complexity of some otherwise  
undecidable or open Diophantine problems. More precisely, 
consider the following: 
\begin{enumerate} 
\renewcommand{\theenumi}{\Roman{enumi}} 
\item{Given a polynomial $f\!\in\!\Z[v,x,y]$, decide the sentence 
$\exists v \; \forall x \; \exists y \; f(v,x,y)\!\stackrel{?}{=}\!0$,\\ 
with all three quantifiers ranging over $\N$ (or $\Z$). } 
\item{Given polynomials $f_1,\ldots,f_m\!\in\!\Z[x_1,\ldots,x_n]$ with 
$m\!\geq\!n$, decide if there\\ is a {\bf rational} solution to 
$f_1\!=\cdots=\!f_m=0$. } 
\end{enumerate} 
We show that problem (I) can be done within $\conp$ for almost all inputs. 
The decidability of problem (I), over $\N$ and $\Z$, was previously 
unknown. We also show that the {\bf Generalized Riemann 
Hypothesis (GRH)} implies that problem (II) can be solved 
within the complexity class $\pp^{\np^\np}$
for almost all inputs, i.e., within the third level 
of the polynomial hierarchy. The decidability 
of problem (II), even in the case $m\!=\!n\!=\!2$, remains open in general.

Along the way, we prove results relating polynomial system 
solving over $\C$, $\Q$, and $\Z/p\Z$. We also prove a result 
on Galois groups associated to sparse polynomial systems which may be of 
independent interest. A practical observation is that the aforementioned 
Diophantine problems should perhaps be avoided in the construction of 
crypto-systems. }  

\begin{article} 

\section{Introduction and Main Results}
\label{sec:intro} 
The negative solution of Hilbert's Tenth Problem 
\cite{undec,hilbert10} has all but dashed earlier hopes of solving 
large polynomial systems over the integers. However, an 
immediate positive consequence is the creation of a 
rich and diverse garden of hard problems with potential applications in 
complexity theory, cryptology, and logic. 
Even more compelling is the question of where the boundary to 
decidability lies. 

From high school algebra we know that detecting and even finding roots in $\Q$ 
(or $\Z$ or $\N$) for polynomials in $\Z[x_1]$ is tractable. (We respectively 
use $\C$, $\R$, $\Q$, $\Z$, and $\N$ for the complex numbers, real numbers, 
rational numbers, integers, and positive integers.) 
However, in \cite{jones9}, James P.\ Jones showed that detecting roots 
in $\N^9$ for polynomials in $\Z[x_1,\ldots,x_9]$ is already 
undecidable.\footnote{This is currently one of the most refined statements 
of the undecidability of Hilbert's Tenth Problem.} 
Put another way, this means that determining the existence of a positive 
integral point on a general algebraic hypersurface of (complex) dimension $8$ 
is undecidable. 

It then comes as quite a shock that decades of number theory 
still haven't settled the complexity of the analogous question for 
algebraic sets  of dimension $1$ through $7$. In fact, even the case 
of plane curves remains a mystery:\footnote{In particular, the major 
``solved'' special cases so far have only extremely ineffective 
complexity and height bounds. (See, 
e.g., the introduction and references of \cite{big}.)} 
As of late 2000, the decidability of 
detecting a root in $\N^2$, $\Z^2$, or even $\Q^2$, for an 
{\bf arbitrary} polynomial in $\Z[x_1,x_2]$, is still completely open.
 
\subsection{Dimensions One and Two} 
To reconsider the complexity of detecting 
integral points on algebraic sets of dimension $\geq\!1$, one can consider 
subtler combinations of quantifiers, and thus subtler questions on the 
disposition of integral roots, to facilitate 
finding  decisive results. For example, Matiyasevich and Julia Robinson have 
shown \cite{matrob,jones81} that sentences 
of the form $\exists u \; \exists v \; \forall x \; \exists y \; 
f(u,v,x,y)\!\stackrel{?}{=}\!0$ (quantified over $\N$), for {\bf arbitrary} 
input $f\!\in\!\Z[u,v,x,y]$, are already undecidable. As another example of the 
richness of Diophantine sentences, Adleman and Manders have shown that 
deciding a very special case of the prefix $\exists\exists$ (quantified 
over $\N$) is $\np$-complete \cite{adleman}: they show $\np$-completeness for 
the  set of $(a,b,c)\!\in\!\N^3$ such that $ax^2+by\!=\!c$ has a solution 
$(x,y)\!\in\!\N^2$.  

However, the decidability of sentences of the form $\exists v \; \forall x 
\; \exists y \; f(v,x,y)\!\stackrel{?}{=}\!0$ (quantified over $\N$ or $\Z$) 
was an open question --- until recently: In \cite{big} it was shown that 
(over $\N$) these sentences can be decided by a Turing machine, 
once the input $f$ is suitably restricted. Roughly speaking, 
deciding the prefix $\exists\forall\exists$ is equivalent to determining 
whether an algebraic surface has a slice (parallel to the 
$(x,y)$-plane) densely peppered with integral 
points.  The ``exceptional'' $f$ not covered by the algorithm of 
\cite{big} form a very slim subset of $\Z[v,x,y]$. 

We will further improve this result by showing that,  
under similarly mild input restrictions, $\exists\forall\exists$ 
can in fact be decided within $\conp$. (This improves a \pspa \, bound 
which appeared earlier in the proceedings version of 
this paper \cite{stoc99}.) To make this more precise, let us write 
any $f\!\in\!\Z[v,x,y]$ as $f(v,x,y)\!=\!\sum c_av^{a_1}x^{a_2}y^{a_3}$, where 
the sum is over certain $a\!:=\!(a_1,a_2,a_3)\!\in\!\Z^3$. 
We then define the 
{\bf Newton polytope of $\pmb{f}$}, $\pmb{\newt(f)}$, as the convex 
hull of\footnote{i.e., smallest convex set in $\R^3$ containing...} 
$\{a \; | \; c_a\!\neq\!0\}$. Also, when we say that a statement 
involving a set of parameters $\{c_1,\ldots,c_N\}$ is true {\bf 
generically}\footnote{ 
We can in fact assert a much stronger condition, but this one suffices for 
our present purposes.}, we will 
mean that for any $M\!\in\!\N$, the statement fails for at most 
$\cO(N(2M+1)^{N-1})$ of the $(c_1,\ldots,c_N)$ lying in $\{-M,\ldots,M\}^N$. 
Finally, for an algorithm with a polynomial $f\!\in\!\Z[v,x,y]$ as input, 
speaking of the {\bf dense encoding} will simply mean 
measuring the input size as $d+\sigma(f)$, where $d$ (resp.\ $\sigma(f)$) is 
the total degree\footnote{i.e., the maximum of the sum of the exponents in any 
monomial term.}  (resp.\ maximum bit-length of a coefficient) of $f$.  
\begin{main}
\label{MAIN:PEPPER} 
Fix the Newton polytope $P$ of a polynomial $f\!\in\!\Z[v,x,y]$  
and suppose that $P$ has at least one integral point in its interior. Assume 
further that we measure input size via the dense encoding. Then, for a generic 
choice of coefficients depending only on $P$, we can decide whether 
$\exists v \; \forall x \; \exists y \; f(v,x,y)\!=\!0$ 
(with all three quantifiers ranging over $\N$ or $\Z$) within $\conp$. 
Furthermore, we can check whether an input $f$ has generic coefficients within 
$\nc$. 
\end{main} 
\noindent 
\begin{rem} 
It is an open question whether membership in 
$\conp$ for the problem above continues to hold relative to the 
{\bf sparse} encoding. 
We will describe the latter encoding shortly. 
Recall also that $\nc\!\subseteq\!\pp\!\subseteq\!\conp$, and 
the properness of each inclusion is unknown \cite{papa}. \qed 
\end{rem} 
The generic choice above is clarified further in section \ref{sec:proof2}. 
It is interesting to note that the exceptional case to our 
algorithm for $\exists\forall\exists$ judiciously contains an extremely hard 
number-theoretic problem: determining the existence of a point in $\N^2$ on an 
algebraic plane curve. (That $\Z[v,y]$ lies in our exceptional locus is easily checked.) 
More to the point, James P.\ Jones has 
conjectured \cite{jones81} that the decidabilities of the prefixes 
$\exists\forall\exists$ and $\exists\exists$, quantified over $\N$, are 
equivalent. Thus, while we have not 
settled Jones' conjecture, we have at least now shown that the decidability of 
$\exists\forall\exists$ hinges on a sub-problem much closer to 
$\exists\exists$. 

It would be of considerable interest to push these techniques further to 
prove a complexity-theoretic reduction from $\exists\forall\exists$ to 
$\exists\exists$, or from $\exists\forall\exists$ to $\forall\exists$. 
This is because these particular reductions would be a first step toward 
reducing $\exists\exists\forall\exists$ to $\exists\exists\exists$, and 
thus finally settling Hilbert's Tenth Problem in three variables. Evidence 
for such a reduction is provided by another result relating (a) the 
size of the {\bf largest} positive integral point on an 
algebraic plane curve with (b) detecting whether an algebraic surface 
possesses {\bf any} integral point: Roughly speaking, it was shown in 
\cite{big} that the computability of the function alluded to in (a) implies 
that the undecidability of $\exists\exists\forall\exists$ occurs only 
in a family of inputs nearly equivalent to $\exists\exists\exists$.  

As for algebraic sets of dimension zero, one can in fact construct \pspa \, 
algorithms to find all {\bf rational} points \cite{stoc99}. 
However, deciding the {\bf existence} of rational points, even for algebraic 
sets of dimension zero, is not yet known to lie within the polynomial 
hierarchy. So let us now consider the latter problem. 

\subsection{Dimension Zero}
\label{sub:zero} 
We will show that deciding feasibility over $\Q$, for most polynomial systems, 
can be done within the polynomial hierarchy, assuming the {\bf 
Generalized\footnote{ The {\bf Riemann Hypothesis (RH)} is an 1859 
conjecture equivalent to a sharp quantitative statement on the 
distribution of primes. GRH can be phrased as a generalization of this 
statement to prime ideals in an arbitrary number field, and further background 
on these RH's can be found in \cite{lago,bs}.} Riemann Hypothesis (GRH)} --- 
a famous conjecture from number theory. To clarify this statement, let us 
first fix some notation and illustrate some of the difficulties presented by  
rational roots of polynomial systems. We will then describe a quantitative 
result depending on GRH before stating our main results on rational roots.  

\noindent 
{\sc Notation} 
{\it  
Let $\mathbf{F}\!:=\!(f_1,\ldots,f_m)$ be a system of polynomials 
in $\Z[x_1,\ldots,x_n]$ and let $\pmb{Z_F}$ be the zero set of $F$ in $\Cn$. 
The {\bf size} of an integer $c$ is 
$\pmb{\size(c)}\!:=\!1+\lceil\log_2(|c|+1)\rceil$. Similarly, 
the {\bf (sparse) size}, $\pmb{\size(F)}$, of the polynomial system $F$ is 
simply the sum of the sizes of all the coefficients {\bf and} exponents in 
its monomial term expansion. \qed }    

To see why it is not entirely trivial to find the rational roots of 
a general $F$ in time polynomial in the sparse size of $F$, consider 
the following two phenomenae: 

\begin{itemize}
\item[{\bf Q$_{\mathbf{1}}$}]{The number of positive integral
roots of $F$ can actually be exponential in $n$: A simple
example is the system $(x^2_1-3x_1+2,\ldots,x^2_n-3x_n+2)$,
with sparse size $\cO(n)$ and root set $\{1,2\}^n$.
Whether the number of rational roots of $F$ can still be exponential
in the sparse size of $F$ for {\bf fixed} $n$ (even $n\!=\!2$!)
is currently unknown. \qed }
\item[{\bf Q$_{\mathbf{2}}$}]{For any {\bf fixed} $n\!>\!1$, the
integral roots of $F$ can have coordinates with bit-length exponential in
$\size(F)$, thus ruling out one possible source of $\np$ certificates:
For example, the system $(x_1-2,x_2-x^d_1,\ldots,x_n-x^d_{n-1})$
has sparse size $\cO(n\log d)$ but has $(2,2^d,\ldots,2^{d^{n-1}})$ as a
root. \qed }
\end{itemize} 
So restricting to deciding the existence of rational
roots, as opposed to explicitly finding them, may be necessary if one wants  
complexity sub-exponential in the sparse size. Indeed, sub-exponential 
bounds are already 
unknown for $m\!=\!n\!=\!2$, and even decidability is unknown in the 
case $F\!:=\!y^2+ax^3+bx+c$ with $a$, $b$, $c$ arbitrary rational 
numbers \cite[ch.\ 8]{sil}, i.e., the case $(m,n)\!=\!(1,2)$. 
So restricting to the case where $Z_F$ is zero-dimensional is 
also crucial. 

On the other hand, when $n\!=\!1$, it is a pleasant surprise that one can 
find {\bf all} rational roots in time polynomial in $\size(F)$ \cite{lenstra}. 
(Note that this is {\bf not} an immediate consequence of the famous 
Lenstra-Lenstra-Lov\'asz factoring algorithm --- the family of examples 
$x^d+ax+b$ already obstructs a trivial application of the latter algorithm.) 
So in order to extend Lenstra's result to general zero-dimensional algebraic 
sets, let us consider an approach {\bf other} than the known 
\pspa \ methods of resultants and  
Gr\"obner bases: reduction modulo specially chosen primes. 

First note that averaging over many primes (as 
opposed to employing a single sufficiently large prime) is essentially 
unavoidable if one wants to use information from reductions modulo primes 
to decide the existence of rational roots.  For example, from basic quadratic 
residue theory 
\cite{hw}, we know that the number of roots $x^2_1+1$ mod $p$ is {\bf not} 
constant for sufficiently large prime $p$. 
Similarly, Galois-theoretic restrictions are also necessary before 
using information mod $p$ to decide feasibility over $\Q$.
\begin{ex}
\label{ex:cool}  
Take $m\!=\!n\!=\!1$ and $F\!=\!f_1\!=(x^2_1-2)(x^2_1-7)(x^2_1-14)$. Clearly, 
$F$ has no rational roots. However, it is easily checked via 
Legendre symbols \cite[ch.\ 9]{apostol} that $F$ has a root 
mod $p$ for {\bf all} primes $p$. In particular, note that the Galois group here does 
not act transitively: there is no automorphism of $\overline{\Q}$ which 
fixes $\Q$ and sends, say, $\sqrt{2}$ to $\sqrt{7}$. \qed 
\end{ex} 

So let us then make the following definition. 
\begin{dfn} 
Let $\pmb{\sigma(F)}$ denote the 
maximum bit-length of any coefficient of the monomial term expansion of 
$F$. Recall that $\pmb{\pi(x)}$ denotes the number of primes $\leq\!x$. 
Let $\pmb{\pi_F(x)}$ be the variation on $\pi(x)$ where we instead count the 
number of primes $p\leq\!x$ such that the mod $p$ reduction of $F$ has a root 
in $\Z/p\Z$. Finally, let $\pmb{N_F(x)}$ be the {\bf weighted} variant of 
$\pi_F(x)$ where we instead count the {\bf total}\footnote{If the number  
of roots in $\Z/p\Z$ of the mod $p$ reduction of $F$ exceeds 
$\delta$, then we add $\delta$ (not $\Omega(p)$) to our total, 
where $\delta$ is as defined in section \ref{sec:wein}.} 
number of distinct roots of the mod $p$ reductions of $F$, summed over all 
primes $p\leq\!x$. \qed 
\end{dfn} 
\noindent 
One can then reasonably guess that behavior 
of the quantities $\frac{\pi_F(x)}{\pi(x)}$ and/or $\frac{N_F(x)}{\pi(x)}$ 
for large $x$ will tell us something about the existence of rational
roots for $F$. 
This is indeed the case, but as we will soon see, the convergence of 
the first quantity to its limit is unfortunately too slow to permit any 
obvious algorithm using sub-exponential work. The second quantity 
will be more important for us algorithmically, so let us give new sharpened 
estimates (depending on GRH) for both quantities.  
\begin{dfn} 
Let $\bO$ and $e_i$ respectively denote the origin and the $i^\thth$ standard 
basis vector of $\Rn$, and normalize $n$-dimensional volume so that the 
standard $n$-simplex (with vertices $\bO,e_1,\ldots,e_n$) has $n$-volume $1$. 
Also let $\#$ denote set cardinality and $V_F\!:=\!\vol_n(Q_F)$, where $Q_F$ is 
the convex hull of the union of $\{\bO,e_1,\ldots,e_n\}$ and the set of all 
exponent vectors of $F$. \qed 
\end{dfn} 
\begin{main} 
\label{MAIN:START} 
Let $K\!:=\!\Q(x_i \; | \; 
(x_1,\ldots,x_n)\!\in\!Z_F \ , \ i\!\in\!\{1,\ldots,n\})$  
and let $r_F$ be\footnote{In \cite{stoc99}, $r_F$ was incorrectly 
defined as the number of rational roots of $F$. } the number of 
maximal ideals in the ring $\Q[x_1,\ldots,x_n]/\langle f_1,\ldots,f_n \rangle$. 
(In particular, $r_F\!\geq\!1$ for $\#Z_F\!\geq\!1$, and for 
$m\!=\!n\!=\!1$ the quantity $r_F$ is just the number of distinct 
irreducible factors of $f_1$ over $\Q[x_1]$.) 
Then the truth of GRH implies the following two statements for all 
$x\!>\!33766$: 
\begin{enumerate}
\item{\ Suppose $\infty\!>\!\#Z_F\!\geq\!2$ and $\mathrm{Gal}(K/\Q)$ acts 
transitively on $Z_F$. Then  
\[ \frac{\pi_F(x)}{\pi(x)}< \left(1-\frac{1}{\#Z_F}
\right)\left(1+\frac{(\#Z_F!+1)\log^2 x + \#Z_F!\cO(\#Z_F\sigma(h_F))\log x} 
{\sqrt{x}}\right)\] }  
\item{\ Suppose $\#Z_F\!\geq\!1$ and $\dim Z_F\!<\!n$. Then independent of 
$\mathrm{Gal}(K/\Q)$, we have 
\[ \frac{\pi_F(x)}{\pi(x)}> \frac{1}{\delta}(r_F-b(F,x)) 
\text{ \ and \ } \left|\frac{N_F(x)}{\pi(x)}-r_F\right|< b(F,x).\] } 
\end{enumerate} 
where $0\!\leq\!b(F,x)\!<\!\frac{4\delta\log^2 x+
\cO(\delta\sigma(\hat{h}_F)(1+n\delta^5/\sqrt{x}))\log
x}{\sqrt{x}}$, 
$0\!\leq\!\sigma(h_F)\!\leq\!\sigma(\hat{h}_F)\!\leq$\\
$\cO(M_F[\sigma(F)+n\log d+\log m])$, $d$ is the maximum 
degree of any $f_i$, $\delta\!\leq\!V_F$, and $M_F$ is no 
larger than the maximum number of lattice points in any translate of 
$(n+1)Q_F$. 
Furthermore, when $m\!\leq\!n$ and $\#Z_F\!<\!\infty$, 
we can replace every occurrence of $\delta$ above with $\#Z_F$. 
Finally, explicit formulae for the asymptotic estimates above 
appear in remarks \ref{rem:ass2} and \ref{rem:ass1} of section 
\ref{sec:wein2}. 
\end{main} 
\begin{rem} 
The polytope volume $V_F$, and even the lattice point count $M_F$, are 
more natural than one might think: $V_F$ is an upper bound on the number of 
irreducible components of $Z_F$ (cf.\ theorem \ref{thm:bezoutkoi} of the next 
section) and $M_F\!=\cO(e^nV_F)$ \cite[sec.\ 6.1.1, 
lem.\ 2 and rem.\ 6]{front}. Furthermore, it is easy to show that 
$V_F\!\leq\!d^n$. In fact, $d^n$ frequently 
exceeds $V_F$ by a factor exponential in $n$ \cite{real,front}. \qed 
\end{rem}
\begin{rem} 
It seems likely that the quantity $\delta$ from theorem \ref{MAIN:START} 
can be replaced by the {\bf affine geometric degree} \cite{cool} and the 
hypotheses $m\!\leq\!n$ and $\#Z_F\!<\!\infty$ dropped. (The affine 
geometric degree agrees with $\#Z_F$ when $\#Z_F\!<\!\infty$ and can be 
significantly less than $V_F$ when $\#Z_F\!=\!\infty$.) This improvement will 
be pursued in future work. \qed
\end{rem} 
 
The upper bound from assertion (1) appears to be new, and the 
first lower bound from
assertion (2) significantly improves earlier
bounds appearing in \cite{hnam,peter} which, when rewritten in the shape
of our bounds, had leading coefficients of $\frac{1}{d^n}$ or worse. 
Also, the special case of the first bound from assertion (2) with 
$m\!\leq\!n$ and 
$F$ forming a reduced regular sequence was independently discovered by Morais 
(see \cite[thm.\ F, pg.\ 11]{morais} or \cite[thm.\ 11, pg.\ 10]{hmps}). 
In this special case, Morais' bound (which depends on the affine geometric 
degree) is asymptotically sharper than our bound when $\#Z_F\!=\!\infty$, and 
our bound is asymptotically sharper when $\#Z_F\!<\!\infty$. We also point out 
that the bounds from \cite[thm.\ F, pg.\ 11]{morais} or \cite[thm.\ 11, pg.\ 
10]{hmps} are stated less explicitly than our formula in remark 
\ref{rem:ass2} of section \ref{sec:wein}, and our 
proof of theorem \ref{MAIN:START} provides a simpler alternative framework 
which avoids the commutative algebra machinery used in \cite{morais,hmps}. 

Part (1) of theorem \ref{MAIN:START} thus presents the main 
difference between ``modular'' feasibility testing over $\C$ and $\Q$: it 
is known \cite[thm.\ 1]{hnam} that the mod $p$ reduction of $F$ has a root 
in $\Z/p\Z$ for a density of primes $p$ which is either positive or zero, 
according as $F$ has a root in $\C$ or not. (See also \cite[sec.\ 2, 
thm.\ 4]{front} for the best current quantitative bound along these lines.) 
The corresponding gap between 
densities is large enough to permit a coarse but fast approximate 
counting algorithm for $\#\mathbf{P}$ to be used to tell the difference, 
thus eventually yielding an \am \, algorithm for feasibility over $\C$ 
recently discovered by Pascal Koiran \cite{hnam}. (We point out that Koiran's 
algorithm also relies on the 
behavior of the function $N_F$, which seems to behave better asymptotically than 
$\pi_F$.) On the other hand, part (1) of theorem \ref{MAIN:START} 
tells us that the mod $p$ reduction of $F$ has a root in $\Z/p\Z$ for 
a density of primes $p$ which is either $1$ or $\leq\!1-\frac{1}{\#Z_F}$ 
(provided $2\!\leq\!\#Z_F\!<\!\infty$), 
and the lower density occurs if $F$ is infeasible over $\Q$ in a 
strong sense. 

Via a $\mathbf{P}^{\np^\np}$ constant-factor 
approximate counting algorithm of Stockmeyer \cite{stock}, we can then 
derive the following result. 
\begin{main}\footnote{This theorem corrects an alleged complexity bound of 
\am, which had an erroneous proof in \cite{stoc99}.} 
\label{MAIN:RIEMANN}
Following the notation and assumptions above, 
assume further that $F$ fails to have a rational root $\Longleftrightarrow 
[Z_F\!=\!\emptyset$ or $\mathrm{Gal}(K/\Q)$  
acts transitively on $Z_F]$. 
Then the truth of GRH implies that deciding whether $Z_F\cap\Qn$ is empty 
can be done within $\pp^{\np^\np}$. Furthermore, we can check the emptiness 
and finiteness of $Z_F$ unconditionally (resp.\ assuming GRH) within \pspa{}  
(resp.\ \am).  
\end{main}

We thus obtain a new arithmetic analogue of Koiran's 
feasibility result over $\C$ \cite{hnam}.  
Indeed, just as we noted for feasibility over $\Q$, the best 
unconditional complexity bound 
for feasibility over $\C$ is \pspa{} \cite{pspace}. However, as we have seen, 
transferring conditional speed-ups from $\C$ to $\Q$ presents some unexpected 
subtleties.

\begin{rem}
The truth of GRH has many other consequences in 
complexity theory. For example, the truth of GRH implies a polynomial time 
algorithm for deciding whether an input integer is prime \cite{miller}, an 
\am\, algorithm 
for deciding whether $Z_F$ is empty \cite{hnam}, and an \am\, 
algorithm for deciding whether $Z_F$ is finite \cite{koiran}. \qed 
\end{rem} 
\begin{rem}
Recall that 
$\mathbf{NP}\!\cup\!\mathbf{BPP}\!\subseteq\!\mathbf{AM}\!\subseteq\!
\corp^{\np}\!\subseteq\!\conp^{\np}\!\subseteq\!\mathbf{P}^{\np^\np}\!
\subseteq\!\cdots\!\subseteq\!\mathbf{PH}\!\subseteq\!
\mathbf{P}^{\#\mathbf{P}}\!
\subseteq\!\mathbf{PSPACE}\!\subseteq\!\mathbf{EXPTIME}$, and the properness 
of each inclusion is unknown \cite{zachos,lab,arith,papa}. \qed 
\end{rem}  
\begin{rem} 
It is quite possible that even without access to an oracle in 
$\np^\np$, the brute-force search implied by the algorithm from 
theorem \ref{MAIN:RIEMANN}, at least for a small number of primes, may 
be more practical than the usual tools of resultants and Gr\"obner bases. 
This remains to be checked extensively. \qed 
\end{rem} 

Let us close with some observations on the strength of our last two theorems: 
First note that our restrictions on the input $F$ are actually rather gentle: 
In particular, if one fixes the monomial term structure of $F$ and 
assumes $m\!\geq\!n$, then it follows easily from the 
theory of resultants \cite{gkz94,introres,gcp} that, for a generic 
choice of the coefficients, $F$ will have only finitely many roots in $\Cn$.  
Furthermore, our hypothesis involving $\mathrm{Gal}(K/\Q)$ holds nearly as 
frequently. 
\begin{thm}
\label{THM:GALOIS} 
Following the notation above, assume $m\!\geq\!n$ and fix the monomial term 
structure of $F$ so that $Z_F\!\neq\!1$ for a generic choice of the 
coefficients. Then, if one restricts to $F$ with integer coefficients 
of absolute value $\leq\!c$, the fraction of such $F$ 
with $\#Z_F\!<\!\infty$ and $\mathrm{Gal}(K/\Q)$ acting transitively  
on $Z_F$ is at least $1-\cO(\frac{\log c}{\sqrt{c}})$. 
Furthermore, we can check whether $\mathrm{Gal}(K/\Q)$ acts 
transitively on $Z_F$ within \expt{} or, if one assumes GRH, within 
$\pp^{\np^\np}$.  
\end{thm}
\noindent
Thus, if $m\!\geq\!n$ and the monomial term structure of $F$ is such that 
$\#Z_F\!\neq\!1$ generically, it 
immediately follows that at least $1-\cO(\frac{\log c}{\sqrt{c}})$ of 
the $F$ specified above have no rational roots. The case 
where the monomial term structure of $F$ is such that $\#Z_F\!=\!1$ 
generically is evidently quite rare, 
and will be addressed in future work. 
\begin{rem} 
A stronger result in the case $m\!=\!n\!=\!1$ (sans complexity 
bounds) was derived by P.\ 
X.\ Gallagher in \cite{gala}. Our more general result above follows from a 
combination of our framework here, the Lenstra-Lenstra-Lovasz (LLL) algorithm 
\cite{lll}, and an effective version of Hilbert's Irreducibility Theorem of 
Stephen D.\ Cohen \cite{cohen}. \qed 
\end{rem} 

Theorems \ref{MAIN:START}--\ref{THM:GALOIS} may thus be of 
independent interest to number theorists, as well as complexity theorists.  
Aside from a geometric trick, the proofs of theorems 
\ref{MAIN:START}--\ref{THM:GALOIS} share a particular tool in common with the 
proof of theorem \ref{MAIN:PEPPER}:  All four proofs make use of some 
incarnation of effective univariate reduction. 

Theorems \ref{MAIN:PEPPER}--\ref{THM:GALOIS} are respectively proved in 
sections 3--6. However, let us first review 
some algorithmic tools that we will borrow from computational algebraic 
geometry and computational number theory. 

\section{Background Tools}
\label{sec:tool} 
We begin with the following elementary fact arising from 
congruences. 
\begin{prop}
\label{prop:dumb}
If $z$ is any {\bf rational} root of 
$\alpha_0+\alpha_1x_1+\cdots+\alpha_dx^d_1\!\in\!\Z[x_1]$, then
$z\!=\!\pm\frac{b}{c}$ for some divisor $b$ of $\alpha_0$ and some divisor $c$
of $\alpha_d$. \qed
\end{prop}

We will also need the following classical fact regarding the factors of a 
multivariate polynomial. 
\begin{lemma} 
\cite[pgs.\ 159--161]{mignotte}
\label{lemma:mignotte}
Suppose $f\!\in\!\Z[t_1,\ldots,t_N]$ has degree $d_i$ with respect
to $t_i$ for all $i$ and coefficients of absolute value $\leq\!c$.
Then $g\!\in\Z[t_1,\ldots,t_N]$ divides $f \Longrightarrow$ the
coefficient of $t^{j_1}_1\cdots t^{j_N}_N$ in $g$ has absolute value $\leq\!
c\prod_i\left(\begin{pmatrix}d_i\\ j_i\end{pmatrix}
\sqrt{(d_i+1)}\right)$, for any
$(j_1,\ldots,j_N)\!\in\![d_1]\times\cdots \times [d_N]$.
In particular, for $N\!=\!1$, $\sigma(g)\!\leq$
\mbox{$\!\sigma(f)+(d_1+\alpha)\log 2$,}
where $\alpha\!:=\!2-\frac{3}{4\log 2}\!<\!0.91798$. \qed
\end{lemma} 

\noindent 
We point out that the last assertion does not appear in \cite{mignotte},
but instead follows easily from Stirling's Estimate \cite[pg.\ 200, ex.\
20]{rudin}. 

We will also need some sufficiently precise quantitative bounds 
on the zero-dimensional part of an algebraic set, e.g., good 
bounds on the number of points and their sizes. A recent bound 
of this type, polynomial in $V_F$, is the following: 
\begin{thm}
\cite[thms.\ 5 and 6]{front} 
\label{thm:bezoutkoi}
Following the notation of section \ref{sub:zero}, there are univariate 
polynomials $P_1,\ldots,P_n,h_F\!\in\!\Z[t]$ with the following properties: 
\begin{enumerate} 
\item{ \ The number of irreducible components of $Z_F$ is bounded above 
by the degree of $h_F$, $\deg h_F$. Furthermore, $\deg P_1,\ldots, 
\deg P_n\!\leq\!  \deg h_F\!\leq\!V_F$, and $\deg h_F\!=\!\#Z_F$ when 
$m\!\leq\!n$ and $\#Z_F\!<\!\infty$.}  
\item{\ $\#Z_F\!<\!\infty \Longrightarrow$ the splitting field of $h_F$ is 
exactly the field 
$K\!=\!\Q[x_i \; \; | \; \; (x_1,\ldots,x_n)\!\in\!\Cn \text{ \ 
is \ a \ root \ of \ } F]$. }  
\item{\ Let $Z'_F$ denote the zero-dimensional part of $Z_F$. 
Then $P_i(x_i)\!=\!0$ for any $(x_1,\ldots,x_n)\!\in\!Z'_F$ 
and any $i\!\in\!\{1,\ldots,n\}$.}   
\item{\ $\sigma(P_1),\ldots,\sigma(P_n)\!\leq\!\sigma(h_F)\!=\!
\cO(M_F[\sigma(F)+n\log d+\log m])$. \qed }  
\end{enumerate} 
\end{thm} 
\begin{rem} 
\label{rem:height} 
Quoting \cite[sec.\ 6.1.1 lem.\ 2 and sec.\ 6.1.3, rem.\ 9]{front},  
we can actually give explicit upper bounds for $\sigma(h_F)$.  
Letting $\mu$ (resp.\ $k$) denote the maximal number of monomial terms 
in any $f_i$ (resp.\ total number of monomial terms in $F$, counting 
repetitions amongst distinct $f_i$), the bounds are as follows:  
\[\mbox{}\hspace{0cm}\log\left\{\frac{16\sqrt{2}}{e^3}
\frac{\sqrt{n+1}}{nV_F}
4^{M_F}\left(n^{3/2}\lceil n \underline{V_F}(\underline{V_F}-1)/4\rceil 
\right)^{V_F} \left(\sqrt{\mu}(c+
\lceil kM_F/2\rceil) \right)^{M_F-V_F} \right\}\]
if $m\!\leq\!n$, or
\[\mbox{}\hspace{-.5cm}\log\left\{\frac{16\sqrt{2}}{e^3}
\frac{\sqrt{n+1}}{nV_F}
4^{M_F}\left(n^{3/2}\lceil n V_F(V_F-1)/4\rceil\right)^{V_F}
\left(\sqrt{\mu}(m\lceil mV_F/2\rceil c+\lceil kM_F/2\rceil
)\right)^{M_F-V_F} \right\} \]
for $m\!>\!n\!\geq\!1$, where $M_F\!\leq\!e^{1/8}
\frac{e^n}{\sqrt{n+1}}V_F+\prod^n_{i=1} (p_i+2)-\prod^n_{i=1} (p_i+1)$, 
and $p_i$ is the length of the projection of $nQ_F$ onto the $x_i$-axis.
(Note that $e^{1/8}\!<\!1.3315$ and $\frac{16\sqrt{2}}{e^3}\!<\!1.127$.) 

Furthermore, if $m\!\leq\!n$ and $\#Z_F\!<\!\infty$, then we can 
replace the underlined occurences of $V_F$ by $\#Z_F$, provided we then add 
an extra summand of $(V_F+\alpha)\log 2$ (with 
$\alpha\!:=\!2-\frac{3}{4\log 2}\!<\!0.91798$) to our bound for $\sigma(h_F)$. 
\qed
\end{rem} 
\begin{rem} 
The true definition of the quantity $M_F$ depends on a particular 
class of algorithms for constructing the {\bf toric resultant} 
(see \cite{front} for further details on $M_F$ and toric resultants). Thus, 
$M_F$ is typically much smaller than the worst-case bound given above. 
\end{rem}  

\noindent 
A preliminary version of the above result was announced in the 
proceedings version of this paper \cite{stoc99}. 
Earlier quantitative results of this type, usually with stronger hypotheses or 
less refined bounds, can be found starting with the work of Joos Heintz 
and his school from the late 80's onward. A good reference for these earlier 
results is \cite{krickpardo} and more recent bounds similar to the one above 
can be found in \cite[prop.\ 2.11]{cool} and \cite[cor.\ 8.2.3]{maillot}. 
There are also more general versions of theorem \ref{thm:bezoutkoi} applying 
even to quantifier elimination over algebraically closed fields, 
but the bounds get looser and the level of generality is greater than we 
need. (These bounds appear in \cite{hnam} and are a corollary of 
results from \cite{fitch}.) 

An immediate corollary of our quantitative result above is the 
following upper bound on $\pi(x)-\pi_F(x)$, which may be of 
independent interest. 
\begin{cor}
\label{cor:easy} 
Following the notation of theorem \ref{thm:bezoutkoi}, assume $F$ has a 
rational root. Then the number of primes $p$ for which the mod $p$ reduction 
of $F$ has {\bf no} roots in $\Z/p\Z$ is no greater than 
$a^*_F\!:=\!n+\sum^n_{i=1}\sigma(P_i)\!=\! 
\cO(nM_F[\sigma(F)+n\log d+\log m])$. 
\end{cor} 

\noindent 
{\bf Proof:} Consider the $i^\thth$ coordinate, $x_i$, of any rational root of 
$F$.  By theorem \ref{thm:bezoutkoi}, and an application of proposition 
\ref{prop:dumb}, the $\log$ of the denominator of $x_i$ (if 
$x_i$ is written in lowest terms) can be no larger than 
$\sigma(P_i)$. In particular, this denominator must have no more 
than $\sigma(P_i)+1$ prime factors, since the only prime power smaller than 
$e$ is $2$. 
Since we are dealing with $n$ coordinates, we can simply sum our last 
bound over $i$ and conclude. \qed 

Let $\pmb{\lie(x)}\!:=\!\int^{x}_{2}\frac{\mathrm{d}t}{\log t}$.  
The following result from analytic number theory will be of 
fundamental importance in our quantitative discussions on prime 
densities. 
\begin{thm} 
\label{thm:lago} 
The truth of RH implies that, for all $x\!>\!2$, 
$\pi(x)$ is within a factor of $1+\frac{7}{\log x}$ 
of $x(\frac{1}{\log x}+\frac{1}{\log^2 x})-\frac{2}{\log 2}$. 
Furthermore, independent of RH, for all $x\!>\!2$, $\lie(x)$ is within a 
factor of $1+\frac{6}{\log x}$ of 
$x(\frac{1}{\log x}+\frac{1}{\log^2 x})-\frac{2}{\log 2}$. \qed 
\end{thm}  
\noindent 
The proof can be sketched as follows: One first 
approximates $\lie(x)$ within a multiple of \mbox{$1+\frac{6}{\log x}$} by 
$x(\frac{1}{\log x}+\frac{1}{\log^2 x})-\frac{2}{\log 2}$, using  
a trick from \cite[pg.\ 80]{apostol}. Then, a (conditional)  
version of the effective Chebotarev Density Theorem, due to Oesterl\'e 
\cite{oyster,bs}, tells us that the truth of RH implies 
\[ |\pi(x)-\lie(x)|<\sqrt{x}\log x, \text{ \ for \ all \ } x\!>\!2. \] 
So, dividing through by 
$x(\frac{1}{\log x}+\frac{1}{\log^2 x})-\frac{2}{\log 2}$ and applying 
the triangle inequality, we obtain our theorem above. 

The remaining facts we need are more specific to the particular 
main theorems to be proved, so these will be mentioned as the need arises. 
\begin{rem} 
Henceforth, we will use a stronger definition of genericity: 
A statement involving a set of parameters $\{c_1,\ldots,c_N\}$ 
holds {\bf generically} iff the statement is true for all 
$(c_1,\ldots,c_N)\!\in\!\C^N$ outside of some {\bf a priori fixed} 
algebraic hypersurface. That this version of genericity implies 
the simplified version mentioned earlier in our theorems is 
immediate from Schwartz' Lemma \cite{schwartz}. \qed 
\end{rem} 

\section{\mbox{Genus Zero Varieties and the Proof of Theorem 
\ref{MAIN:PEPPER}}} 
\label{sec:proof2} 
In what follows, we will make use of some basic algebraic geometry.  
A more precise description of the tools we use can be found in 
\cite{big}. Also, we will always use 
{\bf geometric} (as opposed to arithmetic) genus for  
algebraic varieties \cite{hart}. 

Let us begin by clarifying the genericity condition of theorem 
\ref{MAIN:PEPPER}. Let $Z_f$ be the zero set of $f$. 
What we will initially require of $f$ (in addition to the 
assumptions on its Newton polytope) is that $Z_f$ be  
irreducible, nonsingular, and non-ruled. Later, we will see that a 
weaker and more easily verified condition suffices.  
\begin{rem} 
Ruled surfaces include those surfaces which 
contain an infinite family of lines, for example: planes, cones, 
one-sheeted hyperboloids, and products of a line with a curve. 
More precisely, an algebraic surface $S\!\subseteq\!\Pro^N_\C$ is called {\bf 
ruled} iff there 
is a projective curve $C$, and a morphism $\varphi : S \longrightarrow
C$, such that every fiber of $\varphi$ is isomorphic to $\Pro^1_\C$. 
We then call a surface $S'\!\subseteq\!\C^3$ (the case which concerns us) {\bf 
ruled} iff $S'$ is isomorphic to an open subset of some ruled surface 
in $\Pro^N_\C$. \qed 
\end{rem} 

\begin{lemma}
\label{lemma:geom} 
Following the notation and hypotheses of theorem \ref{MAIN:PEPPER}, 
write $f(v,x,y):=\sum_{(a_1,a_2,a_3)\in A} c_av^{a_1}x^{a_2}y^{a_3}$, 
where $A\cap\{x_i\!=\!0\}\!\neq\!\emptyset$ for all $i$. Then, for a generic 
choice of the coefficients $(c_a)_{a\in A}$, $Z_f$ is irreducible, 
nonsingular, and non-ruled. In particular, for a generic choice of the 
coefficients,  the set $\Sigma_f\!:=\{v_0\!\in\!\C \; | \; \{(x,y)\!\in\!\C^2 
\; | \; f(v_0,x,y)\!=\!0\} 
\text{ \ is \ singular \ or \ reducible}\}$ is finite. 
\end{lemma} 

\noindent 
{\bf Proof:} First note that our hypothesis on $A$ simply prevents the 
coordinate hyperplanes from being subsets of $Z_f$. That $Z_f$ is irreducible 
and nonsingular for a generic choice of coefficients then follows easily from 
the Jacobian criterion for 
singularity \cite{mumford}. (One can even write the conditions explicitly via 
$\cA$-discriminants \cite{gkz94}, but this need not concern us here.)  

That $Z_f$ is also non-ruled generically follows easily from a result of 
Askold G.\ Khovanski relating integral points in Newton polyhedra and 
genera \cite{kho78}: His result, given the hypotheses above, 
implies that $Z_f$ has positive 
genus for a generic choice of the coefficients. (In fact, the only assumptions 
necessary for his result are 
the Newton polytope condition stated in theorem \ref{MAIN:PEPPER} and the 
nonsingularity of $Z_f$.) The classification of algebraic surfaces \cite{beau} 
then tells us that $Z_f$ has positive genus $\Longrightarrow Z_f$ is non-ruled. 

As for the assertion on $\Sigma_f$, assume momentarily that $Z_f$ is 
irreducible, nonsingular, and non-ruled. Then by Sard's theorem 
\cite{hirsch}, $Z_f\cap\{v\!=\!v_0\}$ is irreducible and nonsingular for all 
but finitely many $v_0\!\in\!\C$. Thus, $\Sigma_f$ is finite when $Z_f$ is 
irreducible, nonsingular, and non-ruled.  

Since the intersection of any two open Zariski-dense sets is open and dense, 
we are done. \qed 

\begin{lemma} 
\label{lemma:bound}
Following the notation above, the set of 
$v_0\!\in\!\Z$ such that $\forall x \; \exists y \; 
f(v_0,x,y)\!=\!0$ is contained in $\Sigma_f\cap\Z$, 
whether both quantifiers range over $\N$ or $\Z$. 
Furthermore, $\Sigma_f\cap\N$ finite $\Longrightarrow$ the number of elements 
of $\Sigma_f\cap\Z$, and the size of each such element, is polynomial in the 
dense encoding.  
\end{lemma} 

\noindent 
{\bf Proof:} By Siegel's Theorem 
\cite{wow}, $\forall x\;  \exists y\;  f(v_0,x,y)\!=\!0 
\Longrightarrow Z_f\cap\{v\!=\!v_0\}$ contains a curve of 
genus zero (whether the quantification is over $\N$ or $\Z$).

Now note that for all nonzero $v_0\!\in\!\C$, 
the Newton polytope of $f$ (as a polynomial in {\bf two} variables) is a 
polygon containing an integral point in its interior. So, by Khovanski's 
Theorem \cite{kho78} once again, $Z_f\cap\{v\!=\!v_0\}$ irreducible and 
nonsingular $\Longrightarrow Z_f\cap\{v\!=\!v_0\}$ 
is a curve of positive genus.

Putting together our last two observations, the first part of our lemma 
follows immediately. 

To prove the final assertion, note that the Jacobian 
criterion for singularity \cite{mumford} implies that $\Sigma_f$ is simply 
the set of $v_0$ such that $(v_0,x,y)$ is a complex root of the system of 
equations $(f(v_0,x,y),\frac{\partial f(v_0,x,y)}{\partial x},
\frac{\partial f(v_0,x,y)}{\partial y})$ has a solution $(x,y)\!\in\!\C^2$. 
Thus, $\Sigma_f\cap\N$ finite $\Longrightarrow \Sigma_f$ is a finite set, 
and by theorem \ref{thm:bezoutkoi} we are done. \qed 

Thanks to the following result, we can solve the prefix $\forall\exists$ 
within $\conp$.\\

\noindent 
{\sc Tung's Theorem} 
{\it \cite{tungcomplex} 
Deciding the quantifier prefix $\forall\exists$ (with all quantifiers 
ranging over $\N$ or $\Z$) is $\conp$-complete relative to the dense 
encoding. \qed }\\

\noindent 
The algorithms for $\forall\exists$ alluded in Tung's Theorem are based 
on some very elegant algebraic facts due to James P.\ Jones, Andrzej Schinzel, 
and Shih-Ping Tung. We illustrate one such fact for the case of 
$\forall\exists$ over $\N$. \\

\noindent 
{\sc The JST Theorem} {\it 
\cite{jones81,schinzel,tungcomplex}
Given any $f\!\in\!\Z[x,y]$, we have that
$\forall x\; \exists y\; f(x,y)\!=\!0$
iff all three of the following conditions hold:
\begin{enumerate}
\item{The polynomial $f$ factors into the form
$f_0(x,y)\prod^k_{i=1}(y-f_i(x))$ where $k\!\geq\!1$, $f_0(x,y)\!\in\!\Q[x,y]$
has {\bf no} zeroes in the ring $\Q[x]$, and for all $i$,
$f_i\!\in\!\Q[x]$ and the leading coefficient of $f_i$ is positive.}
\item{$\forall x\!\in\!\{1,\ldots,x_0\} \; \exists
y\!\in\!\N$ such that $f(x,y)\!=\!0$, where $x_0\!=\!\max\{s_1,\ldots,s_k\}$ 
and, for all $i$, $s_i$ is the sum of the squares of the coefficients
of $f_i$.}
\item{Let $\alpha$ be the least positive integer such that
$\alpha f_1,\ldots,\alpha f_k\!\in\!\Z[x]$ and set 
$g_i\!:=\!\alpha f_i$ for all $i$. 
Then the {\bf union} of the solutions of the following $k$ congruences
\begin{eqnarray*}
g_1(x) & \!\equiv & 0 \ \mod \ \alpha \\
 & \vdots & \\
g_k(x) & \!\equiv & 0 \ \mod \ \alpha \text{ \ \ \ \ \ \ \ is 
 \ {\bf all} \ of \ } 
\Z/\alpha\Z. \ \text{\qed} 
\end{eqnarray*} } 
\end{enumerate}
} 

\noindent 
The analogue of the JST Theorem over $\Z$ is essentially the same, save for the 
absence of condition (2), and the removal of the sign check in condition (1) 
\cite{tungcomplex}. 

\medskip
\noindent 
{\bf Proof of Theorem \ref{MAIN:PEPPER}: } Within this proof, we 
will always use the {\bf dense} encoding. Also note that if we 
are quantifying over $\N$, then the roots of $f$ on the coordinate 
hyperplanes can be ignored and we can assume (multiplying by a 
suitable monomial) that the Newton polytope of $f$ intersects 
very coordinate hyperplane. 

Assume $\Sigma_f\cap\N$ is finite. 
{\bf This} will be our genericity hypothesis and 
by lemma \ref{lemma:geom}, and our hypothesis on the Newton polytope of $f$, 
this condition indeed occurs generically. Furthermore, via 
\cite{pspace,neffreif}, we 
can check whether $\Sigma_f$ is finite (and thus whether $\Sigma_f\cap\N$ or 
$\Sigma_f\cap\Z$ is finite) within the 
class $\nc$. It is then clear from lemma \ref{lemma:bound} 
that checking $\exists\forall\exists$ can now be reduced to checking an 
instance of $\forall\exists$ for every $v_0\!\in\!\Sigma_f\cap\N$ (or $v_0\!\in\!\Sigma_f
\cap\Z$). 

Our goal will then be to simply use $\np$ certificates 
for finitely many false $\forall\exists$ sentences, or the 
emptiness of $\Sigma_f\cap\N$ (or $\Sigma_f\cap\Z$), as a single certificate 
of the 
falsity of $\exists\forall\exists$. The emptiness of $\Sigma_f\cap\N$ (or 
$\Sigma_f\cap\Z$) can also be checked within the class $\nc$ \cite{pspace}. So 
by lemma \ref{lemma:bound}, it suffices to assume $\Sigma_f\cap\N$ is 
nonempty and then check that the size of each resulting certificate is 
polynomial in the dense size of $f$.  

Fixing $v_0\!\in\!\Sigma_f\cap\Z$, first note that the dense size of 
$f(v_0,x,y)$ is clearly polynomial in the dense size of $f(v,x,y)$, thanks to 
another application of lemma \ref{lemma:bound}. A certificate of  
$\forall x \; \exists y \; f(v_0,x,y)\!\neq\!0$ (quantified over $\N$) can 
then be 
constructed via the JST Theorem as follows: First, factor $f$ within $\nc$  
(via, say, \cite{bcgw}). 
If $f$ has no linear factor of the form $y-f_i(x)$, then 
we can correctly declare that the instance of $\forall x\; \exists y 
\; f(v_0,x,y)\!\neq\!0$ is true. Otherwise, we attempt to give an 
$x'\!\in\!\{1,\ldots,x_0\}$ such that $f(x',y)$ has 
no positive integral root. Should such an $x'$ exist, lemma 
\ref{lemma:mignotte} tells us that its size will be polynomial in 
$\size(f)$, so $x'$ is an $\np$ certificate. Otherwise, we give a pair 
$(j,t)$ with $1\!\leq\!j\!\leq\!k$ and $t\!\in\!\{0,\ldots,\alpha\}$ such that 
$g_j(t)\!\not\equiv\!0 \ \mod \ \alpha$. Exhibiting such a pair gives a 
negative solution of an instance of the 
{\bf covering congruence} problem, which is known to lie in $\np$ 
\cite{tungcomplex}. 

So we have now proved our main theorem in the case of quantification over 
$\N$. The proof of the case where we quantify over $\Z$ is nearly identical, 
simply using the aforementioned analogue of the JST Theorem over 
$\Z$ instead. \qed 

\begin{rem} 
Note that if $f\!\in\!\Z[v,y]$ then the zero set of $f$ is a ruled 
surface in $\C^3$. From another point of view, the hypothesis of 
theorem \ref{MAIN:PEPPER} is violated since this $P$ has empty interior.  
Deciding $\exists\forall\exists$ for this case then 
reduces to deciding $\exists\exists$, which we've already observed 
is very hard. Nevertheless, Alan Baker has conjectured that the latter 
problem is decidable \cite[sec.\ 5]{jones81}. \qed 
\end{rem} 
\begin{rem} 
The complexity of deciding whether a given surface is ruled 
is an open problem. (Although one can check a slightly weaker condition 
($\#\Sigma_f\!<\!\infty$)
within $\nc$, as noted in our last proof.) It is also interesting to note 
that finding explicit parametrizations of {\bf rational} surfaces (a 
special class of ruled surfaces) appears to be decidable. Evidence 
is provided by an algorithm of Josef Schicho which, while still 
lacking a termination proof, seems to work well in practice \cite{schicho}. 
\qed 
\end{rem} 

\section{Prime Distribution: Proving Theorem \ref{MAIN:START}} 
\label{sec:proof5} 
The proofs of assertions (1) and (2) will implicitly rely on another  
quantitative result on the factorization polynomials, which easily follows 
from Hadamard's inequality \cite{mignotte}. 

\addtocounter{dfn}{1}
{\sc Definition} \arabic{section}.\arabic{dfn}.
{\it Given any polynomial
$f(x_1)\!=\!\alpha_0+\alpha_1x_1+\cdots+\alpha_Dx^D_1$, 
we define:}\\
\vbox{
\begin{wrapfigure}{l}{3in}
\scalebox{.8}[1]{
\mbox{$\pmb{\Delta_f}:=\frac{(-1)^{D(D-1)/2}}{\alpha_D}$
\begin{footnotesize}
$\mathrm{DET}
\begin{bmatrix}
\alpha_0 & \cdots & \alpha_D & 0 & \cdots & 0 & 0 \\
0     & \alpha_0 & \cdots & \alpha_D & 0 & \cdots & 0\\
\vdots &  \ddots  & \ddots &    &  \ddots  & \ddots  & \vdots   \\
0      & \cdots & 0 & \alpha_0 & \cdots & \alpha_D & 0 \\
0      & 0  & \cdots & 0 & \alpha_0 & \cdots & \alpha_D \\
\alpha_1 & \cdots & D\alpha_D & 0 & \cdots & 0 & 0 \\
0     & \alpha_1 & \cdots & D\alpha_D & 0 & \cdots & 0\\
\vdots &  \ddots  & \ddots &    &  \ddots  & \ddots  & \vdots   \\
0      & \cdots & 0 & \alpha_1 & \cdots & D\alpha_D & 0 \\
0      & 0  & \cdots & 0 & \alpha_1 & \cdots & D\alpha_D
\end{bmatrix}$
\end{footnotesize}},}  
\end{wrapfigure}
\mbox{}\\
{\it ...where the first $D-1$ (resp.\ last $D$) rows of the matrix correspond
to the coefficients of $f$ (resp.\ the derivative of $f$).
The quantity $\Delta_f$ is also known as the {\bf discriminant of}
$\mathbf{f}$, and vanishes only for polynomials with repeated roots 
\cite{gkz94}. \qed }
}
\mbox{}\\
\vspace{-.5cm}
\mbox{}\\
\begin{lemma}
\label{lemma:disc}  
Suppose $g\!\in\!\Z[x_1]$ is square-free and $\delta\!:=\!\deg g$. Then
\[\log |\Delta_g|\!\leq\!(2\delta-1)\sigma(f)+ 
\frac{2\delta-1}{2}\log(\delta+1)+\frac{\delta}{2}\log(\delta(2\delta+1)/6). 
\text{ \ \qed}\]
\end{lemma} 

The last and most intricate result we will need is the 
following refined effective version of the primitive element theorem.  
\begin{thm} 
\cite[thm.\ 7]{front} 
\label{thm:prim}
Following the notation of theorem \ref{thm:bezoutkoi}, 
one can pick $\hat{h}_F\!\in\!\Z[t]$ (satisfying all the properties of $h_F$ 
from theorem \ref{thm:bezoutkoi}), so that there also 
exist $a_1,\ldots,a_n\!\in\!\N$ and $h_1,\ldots,h_n\!\in\!\Z[t]$ with the
following properties: 
\begin{enumerate}
\item{\ The degrees of $h_1,\ldots,h_n$ are all bounded above
by $\deg(\hat{h}_F)\!\leq\!V_F$.}
\item{\ For any root $(\zeta_1,\ldots,\zeta_n)\!\in\!Z'_F$ of $F$, there 
is a root $\theta$ of $\hat{h}_F$ such that  
$\frac{\ h_i(\theta)}{a_i}\!=\!\zeta_i$ for all $i$. }
\item{\ For all $i$, both $\log a_i$ and
$\sigma(h_i)$ are bounded above by $\cO(V^5_F\sigma(\hat{h}_F))$ and 
$\sigma(\hat{h}_F)\!=\!\cO(\sigma(h_F))$. \qed } 
\end{enumerate}
\end{thm}
\begin{rem} 
\label{rem:loga} 
Quoting \cite[sec.\ 6.1.5, rem.\ 11]{front},  
we can actually make the asymptotic bounds above completely explicit:  
\[ \sigma(h_i)\!\leq\!(2\delta^2-2\delta+1)\sigma(r)+(2\delta^2+1)
\sigma(\hat{h}_F) 
+\log[(\delta^2+1)^{\delta^2}(\delta+1)^{\delta+1}(\delta^2-\delta+1)] \]
and 
\[ \log a_i\!\leq\!\delta(\delta-1)\sigma(r)+(\delta^2+1)\sigma(\hat{h}_F)
+\frac{1}{2}\log[(\delta^2+1)^{\delta^2}(\delta+1)^\delta],\] 
where $\sigma(r)\!\leq\!\frac{\delta^2-\delta+2}{2}\log(B^2_1
+\delta(\delta-1)/2)$, 
$B_1\!:=\!\left( \frac{4\cdot {16}^{\delta+1}}{e^{9/4}}\cdot
\sqrt{(\delta+1)^5} \right)^{\delta-1} e^{2(\delta-1)\sigma(\hat{h}_F)}$, 
$\delta\!:=\!\max \deg h_i\!\leq\!\deg \hat{h}_F\!\leq\!V_F$, 
$\sigma(\hat{h}_F)\!\leq\!\sigma(h_F)+\delta'\log(2n+1)+(V_F+\alpha)\log 2$, 
$\delta'\!\leq\!V_F$, $\sigma(h_F)$ is bounded above as in remark 
\ref{rem:height} of section \ref{sec:tool}, and 
$\alpha\!:=\!2-\frac{3}{4\log 2}\!<\!0.91798$. 
(So $\log a_i$ actually admits an upper bound about half as 
large as the bound for $\sigma(h_i)$.) 

Furthermore, when $m\!\leq\!n$ and $\#Z_F\!<\!\infty$, 
we can replace every occurence of $\delta$ and $\delta'$ above by $\#Z_F$. \qed 
\end{rem} 
\begin{rem} 
Earlier quantitative results of this type, e.g., those applied 
in \cite{hnam}, had looser and less explicit bounds which were polynomial in 
$d^{n^{\cO(1)}}$.  
\qed 
\end{rem} 

\subsection{Proving Assertion (2) of Theorem \ref{MAIN:START}}\mbox{}\\
\label{sec:wein} 
First let us recall the following refined version of an important result due 
to Weinberger.  
\begin{thm} 
\label{thm:wein} 
Following the notation of lemma \ref{lemma:disc}, suppose 
$g\!\in\!\Z[x_1]$ has degree $\delta$ and {\bf no} factors of multiplicity 
$>\!1$. Then the truth of GRH implies that 
\[ \left|\frac{N_g(x)}{\pi(x)}-r_g\right|< 
\frac{2\sqrt{x}\log(|\cD_g|x^\delta)+\delta\log|\cD_g|} 
{\lie(x)}, \text{ \ for \ all \ } x\!>\!2. \text{ \ \qed } \] 
\end{thm} 
\noindent 
The original version from \cite{weinberger} had an unspecified constant in 
place of the $2$. The version above follows immediately from Weinberger's 
original proof, simply using a stronger version of effective Chebotarev than 
he used, i.e., one replaces theorem 1.1 of \cite{lago} by a 
result of Oesterl\'e \cite{oyster} (see also theorem 8.8.22 of \cite{bs}). 

The second (harder) bound of assertion (2) of Theorem \ref{MAIN:START} 
is then just a simple corollary of theorems \ref{thm:bezoutkoi} and 
\ref{thm:wein}. The first bound is an even simpler corollary of the second 
bound. 
 
\noindent 
{\bf Proof of Assertion (2):} 
By theorems \ref{thm:bezoutkoi} and \ref{thm:prim}, it 
immediately follows that $r_F\!=\!r_{\hat{h}_F}$. 
(Note that $\hat{h}_F$ is square-free by construction.)  
It also follows easily that the mod $p$ reduction of $F$ has a root in $\Z/p\Z 
\Longrightarrow$ the mod $p$ reduction of $\hat{h}_F$ has a root in $\Z/p\Z$. 
Furthermore, theorem \ref{thm:prim} tells us that a sufficient 
condition for the converse assertion is that $p$ not divide any of the $a_i$ 
(the denominators in our rational univariate representation of $Z_F$). 
We thus obtain $0\!\leq\!N_{\hat{h}_F}(x)-
N_F(x)\!\leq\!\delta\sum^n_{i=1}(\log a_i+1)$, 
for all $x\!>\!0$, where $\delta\!:=\!\deg \hat{h}_F$. 

Assume henceforth that $x\!>\!2$. We then have  
\[ \left|\frac{N_F(x)}{\pi(x)}-r_F\right|\!\leq\!\left|\frac{N_{\hat{h}_F}(x)}
{\pi(x)}-r_{\hat{h}_F}\right|+ 
\frac{\delta(\sum^n_{i=1}\log a_i + n)}{\pi(x)}.\] 
Combining theorem \ref{thm:wein} and 
Oesterl\'e's conditional bound on $|\pi(x)-\lie(x)|$, we thus obtain 
that the truth of GRH implies 
\[\left|\frac{N_F(x)}{\pi(x)}- r_F\right|\!<\!\frac{2\sqrt{x}
\log(|\cD_{\hat{h}_F}| 
x^{\delta})+\delta\log|\cD_{\hat{h}_F}|}
{\lie(x)}+\left(1+\frac{\sqrt{x}\log x}{\lie(x)}
\right)\frac{ \delta(\sum^n_{i=1}\log a_i + n)}{\lie(x)}.\] 

By theorem \ref{thm:lago}, and the fact that 
$\frac{(\log^3 x)(1+6/\!\log x)}{\sqrt{x}(\log x+1)-\frac{2}{\log 2}\log^2 x}
\!<\!1$ for all $x\!>\!33766$, we then obtain 
\[\left|\frac{N_F(x)}{\pi(x)}- r_F\right|\!<\!\frac{2\sqrt{x}
\log(|\cD_{\hat{h}_F}| 
x^\delta)+\delta\log|\cD_{\hat{h}_F}|+2\delta(\sum^n_{i=1}\log a_i + n)}
{\lie(x)},\] 
for all $x\!>\!33766$. The second bound from assertion (2) then follows 
immediately from lemma \ref{lemma:disc}, theorem \ref{thm:bezoutkoi}, and the 
fact that $\frac{\lie(x)}{x/\!\log x}\!<\!(1+4/\!\log x)^2$ (applying  
theorem \ref{thm:lago} one last time). 

The first bound of assertion (2) follows immediately from the 
second bound via a simple application of the triangle inequality and 
the inequality $N_F(x)\!\leq\!\delta\pi_F(x)$. \qed 
\begin{rem} 
\label{rem:ass2} 
Carrying out the last step in detail (and observing that 
$(1+4/\!\log x)^2\!<\!2$ for all $x\!>\!33766$) it is clear that 
the asymptotic bound on $b(F,x)$ can be replaced 
by the following explicit quantity: 
\[\frac{4\delta\log^2 x+
\left(4\log|\cD_{\hat{h}_F}|+\frac{2\delta\left(\log |\cD_{\hat{h}_F}|+2n
+2\sum^n_{i=1}\log a_i\right)}
{\sqrt{x}}\right)\log x}{\sqrt{x}},\]
where $\log |\Delta_{\hat{h}_F}|\!\leq\!(2\delta-1)\sigma(\hat{h}_F)
+\frac{2\delta-1}{2}\log(\delta+1)+
\frac{\delta}{2}\log(\delta(2\delta+1)/6)$, 
$\delta\!:=\!\deg \hat{h}_F\!\leq\!V_F$, and  
$\hat{h}_F$ and $\log a_i$ are 
as in theorem \ref{thm:prim} and remark \ref{rem:loga} of section 
\ref{sec:proof5}. 

Furthermore, via \cite[sec.\ 6.1]{front}, we can conclude that every 
occurence of $\delta$ can be replaced by $\#Z_F$ when  
$m\!\leq\!n$ and $\#Z_F\!<\!\infty$. \qed  
\end{rem} 

\subsection{Proving Assertion (1) of Theorem \ref{MAIN:START}}\mbox{}\\ 
\label{sec:wein2} 
Here we will need the following result dealing with the density of 
primes for which the mod $p$ reduction of $F$ has a root in $\Z/p\Z$. 
This theorem may be of independent interest to computational number theorists. 
\begin{thm}
\label{thm:salvation}
Following the notation of theorem \ref{MAIN:START}, assume $\#Z_F\!<\!\infty$ 
and let $j_F$ 
be the fraction of elements of $\mathrm{Gal}(K/\Q)$ which fix 
at least one root of $F$. Then the truth of GRH implies that 
\[\left|\frac{\pi_F(x)}{\pi(x)}- j_F\right|\!<\!\frac{j_F(V_F!+1)\log^2 x 
+2\left(j_FV_F!\log|\cD_g| + \frac{\sigma(h_F)+1}{\sqrt{x}}\right)\log x}
{\sqrt{x}},\] 
for all $x\!>\!33766$, where $h_F$ is the polynomial from 
theorem \ref{thm:bezoutkoi} and $g$ is the square-free part of $h_F$. 
\end{thm} 

\noindent 
{\bf Proof:} 
Let $g$ be the square-free part of the polynomial $h_F$ from 
theorem \ref{thm:bezoutkoi} and let $j_g$ be the fraction of elements of the 
Galois group of $g$ (over $\Q$) 
which fix at least one root of $g$, where $g$ is the square-free part of 
the polynomial $h_F$ from theorem \ref{thm:bezoutkoi}. By 
essentially the same argument as the beginning of the proof of assertion (1), 
we obtain $j_F\!=\!j_g$. Similarly, we also obtain 
$0\!\leq\!\pi_g(x)-\pi_F(x)\!\leq\!\sigma(h_F)+1$ for all 
$x\!>\!2$.

Note that $j_g$ is also the fraction of elements of the Galois group which 
give permutations (of the roots of $g$) possessing a fixed point. 
Oesterl\'e's (conditional) version 
of effective Chebotarev \cite{oyster,bs} then tells us\footnote{ 
His result is actually stated in 
terms of conjugacy classes, but since the number of fixed points of a 
Galois group element is stable under conjugacy, we can simply sum over 
conjugacy classes.} that the truth of GRH implies 
$\left| \pi_g(x)-j_g\lie(x)\right|\leq j_g\sqrt{x}(2\log|\cD|+{\cd}\log x)$, 
where $\cD$ is the discriminant of the splitting field of $g$ and 
${\cd}$ is the degree of this field extension over $\Q$. Letting 
$\delta\!:=\!\deg g$ (which is exactly $\#Z_F$ by construction), basic Galois 
theory tells us that $\cd\!\leq\!\#Z_F!$. 

By Oesterl\'e's conditional bound on $|\pi(x)-\lie(x)|$ we then  
obtain 
\[ \left| \pi_g(x)-j_g\pi(x)\right|\leq j_g\sqrt{x}(2\log|\cD|+(\cd+1)\log x).
\] 
Following essentially the same reasoning as the proof 
of assertion (2) we then obtain 
\[\left|\frac{\pi_F(x)}{\pi(x)}- j_F\right|\!<\!\frac{j_g(\cd+1)\log^2 x 
+2\left(j_g\log|\cD| +\frac{\sigma(h_F)+1}{\sqrt{x}}\right)\log x}{\sqrt{x}},\] 
for all $x\!>\!33766$. Using the fact that $|\cD|\!\leq\!|\cD_g|^{\cd}$ 
\cite[pg.\ 259]{bs}, and applying lemma 
\ref{lemma:disc}, we are done. \qed 

Of course, we must now estimate the quantity $j_F$. Fortunately, 
a good upper bound has already been derived by Peter J.\ Cameron 
and Arjeh M.\ Cohen, in answer to a 1991 question of Hendrik 
W.\ Lenstra. 
\begin{thm}
\label{thm:group} 
Suppose $G$ is any group acting transitively and faithfully 
on a set of $N$ elements and $j_G$ is the fraction of elements of $G$ 
with at least one fixed-point. Then $j_G\!\leq\!1-\frac{1}{N}$. \qed 
\end{thm}
\noindent 
The proof occupies the second page of \cite{camcoh} and 
requires only some basic group representation theory.\footnote{ Their 
paper actually dealt with finding a {\bf lower} bound for the quantity 
$1-j_G$. } The upper bound is tight, but completely classifying the next 
lower values of $j_G$ currently requires the classification of finite simple 
groups \cite{guralwan}. The latter classification will {\bf not} be necessary 
for our results.  

\medskip 
\noindent 
{\bf Proof of Assertion (1):} Following the notation of our last proof, 
recall that $g$ is the square-free part of the polynomial $h_F$ from theorem 
\ref{thm:bezoutkoi}. Then by assumption, $V_F\!\geq\!\#Z_F\!\geq\!2$ 
and $\delta\!=\!\#Z_F$. Furthermore, by theorems \ref{thm:bezoutkoi} 
and \ref{thm:group}, $j_F\!\leq\!1-\frac{1}{\#Z_F}$. So by 
theorem \ref{thm:salvation} we are done. \qed 
\begin{rem} 
\label{rem:ass1} 
From our proofs above we easily see that the asymptotic bound from assertion 
(1) can be replaced by the following explicit quantity: 
\[ \left(1-\frac{1}{\#Z_F}
\right)\left(1+\frac{(\#Z_F!+1)\log^2 x + 2\left(\#Z_F!\log|\cD_g|+
\frac{\#Z_F}{\#Z_F-1}\cdot \frac{\sigma(h_F)+1}{\sqrt{x}}\right)\log x} 
{\sqrt{x}}\right),\] 
where $g$ is as in our proof above, 
$\log|\cD_g|\!\leq\!2(\delta-1)(\sigma(h_F)+(V_F+\alpha)\log 2)+
\frac{2\delta-1}{2}\log(\delta+1)+\frac{\delta}{2}\log(\delta(2\delta+1)/6)$ 
(thanks to lemmata \ref{lemma:mignotte} and \ref{lemma:disc}), 
$\alpha\!:=\!2-\frac{3}{4\log 2}\!<\!0.91798$, 
and $\sigma(h_F)$ is bounded as in remark \ref{rem:height} of section 
\ref{sec:tool}. \qed  
\end{rem} 

\section{The Proof of Theorem \ref{MAIN:RIEMANN}} 
\label{sec:proof1}
Our algorithm essentially boils down to checking whether $r_F\!\geq\!2$ 
or $r_F\!=\!1$, following the notation of theorem \ref{MAIN:START}. Via our 
initial assumptions on $F$, we will see that this is the same as checking 
whether $F$ as a rational root or not. 
\begin{rem} 
It is at this point that we must slightly alter our defintion 
of $N_F$: As we sum the number of roots in $\Z/p\Z$ of the mod $p$ 
reductions of $F$, we instead add $V_F$ to our total for each $p$ where this 
number of roots exceeds $V_F$. This ensures that $N_F$ can actually be 
computed within $\#\pp$, since $V_F$ can be computed within $\#\pp$ (see 
below).  It is unknown whether the same is true for the quantity $\delta$ in 
our initial definition of $N_F$. \qed 
\end{rem} 

Our algorithm proceeds as follows:  
First check whether $Z_F$ is empty. If so, then we immediately 
know that $Z_F\cap\Q^n$ is empty and we are done. 
Otherwise, approximate $N_F(M)$ and $\pi(M)$ within a factor of $\frac{9}{8}$, 
where $M$ is an integer sufficiently larger than $33766$ so that 
$b(F,M)\!<\!\frac{1}{10}$. Respectively calling these approximations 
$\bar{N}$ and $\bar{\pi}$, we then do the following: 
If $\bar{N}\!\leq\!(\frac{9}{8})^2\bar{\pi}$, declare $Z_F\cap\Q^n$ 
empty. Otherwise, declare $Z_F\cap\Q^n$ nonempty. 

That our algorithm works is easily checked. First note that  
$\bar{N}\!\leq\!(\frac{9}{8})^2\bar{\pi} \Longleftrightarrow 
\frac{N_F(M)}{\pi(M)}\!\leq\!(\frac{9}{8})^4$. So by theorem \ref{MAIN:START},  
our assumption on $b(F,M)$ implies that the last inequality occurs iff 
$r_F\!=\!1$. (Note that we need GRH at this point.) 
Via theorem \ref{thm:prim}, and our earlier proofs, we know that  
$r_F\!=\!r_{\hat{h}_F}$. So by \cite[thm.\ 4.14]{jacob}, we have 
that $\mathrm{Gal}(K/\Q)$ acts transitively on $Z_F$ iff $\hat{h}_F$ is 
irreducible over 
$\Q$ (or equivalently, $r_F\!=\!r_{\hat{h}_F}\!=\!1$). So by our initial 
assumptions on 
$F$, $r_F\!=\!1$ iff $F$ has no rational roots. Thus, we now need only check 
the complexity of our algorithm. 

That the emptiness and finiteness of $Z_F$ can be checked within 
\pspa{} unconditionally goes back to \cite{pspace}. That the truth of GRH 
implies both bounds can be lowered to \am{} is proved respectively in 
\cite{hnam} and \cite{koiran}. So now we need only check the 
complexity of computing $M$, $\bar{N}$, and $\bar{\pi}$. 

It follows immediately from \cite{pratt} that $N_F(x)$ 
and $\pi(x)$ can be computed within $\#\pp$. Also, 
via \cite{volcomplex}, $V_F$ can be computed within $\#\pp$ as well. 
Furthermore, via theorems \ref{MAIN:START} and \ref{thm:bezoutkoi} 
(and the fact that $0\!\leq\!\log V_F\!\leq\!n\log d$), the number 
of bits of $M$ is polynomial in the size of $F$. So by \cite{stock}, 
$M$, $\bar{N}$, and $\bar{\pi}$ can be computed within $\pp^{\np^\np}$. 
Therefore, our algorithm runs within $\pp^{\np^\np}$, assuming GRH. \qed 
\begin{rem} 
It is an open problem whether theorem \ref{MAIN:RIEMANN} continues to 
hold under the weaker condition that the {\bf real} dimension of $Z_F$ is at 
most zero. \qed 
\end{rem} 

\section{The Proof of Theorem \ref{thm:galois}} 
\label{sec:final} 
If $m\!>\!n$ then it follows easily from Schwartz' Lemma \cite{schwartz} 
that $F$ has {\bf no} roots for at least a fraction of 
$1-\cO(\frac{1}{c})$ of our $F$. 
So we can assume $m\!=\!n$. 

Consider now the toric resultant, $\cR$, of 
$f_1,\ldots,f_n$ and $u_0+u_1x_1+\cdots+u_nx_n$. 
(The classical resultant of Macaulay would suffice to 
prove a weaker version of our theorem here for a more  
limited family of monomial term structures.) Then, for indeterminate 
coefficients, $\cR$ is a nonzero irreducible polynomial over $\Z$ adjoin $u_0,
\ldots,u_n$ and the coefficients of $F$. More importantly, if the 
coefficients of $F$ are {\bf constants}, $\cR$ is divisible by 
$u_0-(\zeta_1u_1+\cdots+\zeta_nu_n)$, for any root 
$(\zeta_1,\ldots,\zeta_n)\!\in\!\C^n$ of $F$. 

If it happens that $\cR$ (in fully symbolic form) is the constant $1$, then 
it follows from the degree formula for the toric resultant \cite{gkz94} 
that $Z_F$ is empty for a generic choice of the coefficients and 
there is nothing to prove. 
So let us assume $\cR$ is not identically $1$ in its full symbolic form. 

By \cite{cohen} it then follows that a fraction of at most 
$\cO(\frac{\log c}{\sqrt{c}})$ of the $F$ whose coefficients 
are {\bf rational} numbers of (absolute multiplicative) height $\leq\!c$ 
result in $\cR$ being a 
reducible polynomial over $\Q[u_0,\ldots,u_n]$. By rescaling, 
this easily implies that at most $\cO(\frac{\log c}{\sqrt{c}})$ of 
the $F$ whose coefficients are {\bf integers} of absolute value $\leq\!c$ 
result in $\cR$ being reducible over $\Q[u_0,\ldots,u_n]$. 

We now observe (say from 
\cite[sec.\ 6]{front}) that the polynomial $h_F$ from theorem 
\ref{thm:bezoutkoi} is nothing more than the resultant $\cR$, for suitably 
chosen $u_1,\ldots,u_n$. (So in particular, $\cR$ irreducible and nonzero 
$\Longrightarrow \#Z_F\!<\!\infty$.) So let us apply the 
Effective Hilbert Irreducibility Theorem from \cite{cohen} one more 
time to obtain such a choice of $u_1,\ldots,u_n$. 

We then obtain that the fraction of our $F$ for which 
$\#Z_F\!<\!\infty$ and $h_F$ is irreducible over $\Q$ is at least 
$1-\cO(\frac{\log c}{\sqrt{c}})$. By \cite[thm.\ 4.14]{jacob}, $h_F$ is 
irreducible iff its Galois group acts transitively on its roots. So by theorem 
\ref{thm:bezoutkoi}, our first assertion is proved. 

That $\mathrm{Gal}(K/\Q)$ acts transitively on $Z_F$ can be checked 
within $\pp^{\np^\np}$ (assuming GRH) is already clear from the proof 
of theorem \ref{MAIN:RIEMANN}. To obtain the unconditional 
complexity bound, it clearly suffices to factor $h_F$ within \expt{}  
and see whether $h_F$ is irreducible. Since theorem \ref{thm:bezoutkoi} 
tells us that the dense size of $h_F$ is exponential in $\size(F)$, we can 
conclude via an application of the polynomial-time LLL factoring algorithm 
from \cite{lll}. \qed 

\section*{Acknowledgements} 
The author thanks Felipe Cucker, Jan Denef, Michael Fried, Teresa Krick, 
Jeff Lagarias, Luis-Miguel Pardo, and Bjorn Poonen for some very useful 
discussions, in person and via e-mail. 
In particular, Jan Denef pointed out the excellent 
reference \cite{friedjarden}, and Michael Fried helped confirm a 
group-theoretic hope of the author (theorem \ref{thm:group}). Special thanks 
go to Pascal Koiran for pointing out errors in earlier versions of theorems  
\ref{MAIN:START} and \ref{MAIN:RIEMANN}. 

This paper is dedicated to Gretchen Davis, a remarkable educator who 
first inspired the author's interest in mathematics.

\end{article} 
\end{document}